\documentclass[11pt,reqno]{amsart}
\usepackage[top=90pt,bottom=90pt,left=90pt,right=90pt]{geometry}
\usepackage{bbm,amsmath,amssymb,amsfonts,graphicx,color,subfig,mathtools,verbatim,MnSymbol}
\usepackage{hyperref,float}
\usepackage[ruled,lined,linesnumbered]{algorithm2e}
\usepackage[dvipsnames]{xcolor}
\usepackage{comment}
\usepackage{graphicx}
\usepackage[shortlabels]{enumitem}
\usepackage{tikz}
\usepackage{float}
\usetikzlibrary {arrows.meta,bending} 
\hypersetup{
}

\usepackage[
maxnames=99,
backend=biber,
style=numeric,
sorting=nyt,
giveninits=true
]{biblatex}
\numberwithin{equation}{section}
\usepackage{aliascnt}
\newtheorem{theorem}{Theorem}[section]
\newtheorem{lemma}[theorem]{Lemma}
\newtheorem{remark}[theorem]{Remark}

\newtheorem{corollary}[theorem]{Corollary}
\newtheorem{proposition}[theorem]{Proposition}
\newtheorem{definition}[theorem]{Definition}
\newtheorem{example}[theorem]{Example}

\def\NN{\mathbb{N}}
\def\ZZ{\mathbb{Z}}
\def\QQ{\mathbb{Q}}
\def\RR{\mathbb{R}}

\def\N{\mathcal{N}}
\def\M{\mathcal{M}}
\def\S{\mathcal{S}}
\def\C{\mathcal{C}}
\def\R{\mathcal{R}}

\newcommand{\setn}[1]{\left\{#1\right\}}

\newcommand{\setcond}[2]{\left\{#1 \,\colon\, #2\right\}}

\newcommand{\suchthat}{\;\ifnum\currentgrouptype=16 \middle\fi|\;}


\begin{document}
\title[Symmetric Chemical Reaction Networks: Realizations and Steady States]{Symmetric Chemical Reaction Networks: \\ Realizations and Steady States}
\author{Sebastian Debus}
\address{Universität Konstanz, Fachbereich Mathematik und Statistik, 78457 Konstanz, Germany}
\email{sebastian.debus@uni-konstanz.de}
\author{Kurt Klement Gottwald}
\address{Technische Universität Chemnitz, Fakultät für Mathematik, 09107 Chemnitz, Germany}
\email{kurt-klement.gottwald@mathematik.tu-chemnitz.de}
\author{Nils Gubela}
\address{Department of Mathematics and Computer Science, Freie Universität Berlin, Berlin, Germany, International Max-Planck Research School “Biology and Computation” (IMPRS-BAC), Max-Planck Institute for Molecular Genetics, Berlin, Germany }
\email{nils.gubela@fu-berlin.de}

\begin{abstract}
We investigate symmetric chemical reaction networks and their steady states. 
From classical equivariant dynamical systems and network theory it is known that the mass-action ODE of any symmetric chemical reaction network is equivariant. For symmetries that permute the species, we prove that any equivariant mass-action ODE is realized by a symmetric chemical reaction network and construct a sparsest such realization. We apply the half-degree principle to show that a symmetric CRN that has a steady state also has one with few distinct coordinates.
\end{abstract}
\maketitle

\section{Introduction}
A \emph{chemical reaction network} (CRN) $\mathcal{N}$ is a triple $(\mathcal{S},\mathcal{C},\mathcal{R})$ consisting of a finite set of pairwise distinct \emph{species} $\mathcal{S}=\{X_1,\ldots,X_n\}$, a finite set of \emph{complexes} $\mathcal{C}$, i.e., a subset of $\setcond{a_1X_1+\ldots + a_n X_n }{ a_1,\ldots,a_n \in \NN_0}$, and a finite set of \emph{reactions} $\mathcal{R}$, i.e., a subset of $\setcond{A \to B}{A, B \in \mathcal{C}, A \neq B}$. Thus, a CRN can be thought of as a directed graph with vertices $\mathcal{C}$.
The edges of this graph represent the pathways of chemical transformations, where each reaction is typically associated with a rate function that governs the system's temporal evolution. In the most common framework, known as \emph{mass-action kinetics} \cite{feinberg2019foundations}, the network topology is translated into a system of autonomous ordinary differential equations describing the concentrations of the species over time. 

One of the central goals of CRN theory is to characterize the connection between network structure and the dynamics it induces, particularly at steady state.
This requires finding real roots in the non-negative orthant of a system of nonlinear multivariate polynomials. As the size of the network grows, standard symbolic computation techniques for solving these systems, such as Gröbner bases, often become computationally intractable due to algorithmic complexity \cite{arXiv:2507.23728}.

Many chemical reaction networks possess inherent symmetries that can be exploited to reduce the algorithmic complexity. We distinguish between two types of symmetries: graph theoretical symmetries of the CRN (i.e. permuting vertices and edges), and symmetries acting directly on the species. 
It is widely known that symmetric networks induce equivariant dynamical systems \cite{golubitsky2003symmetry}. 
Symbolic algorithms have exploited these symmetries via representation theory, invariant rings, or the decomposition of the ODE into equivariants \cite{chossat2000methods, Gatermann00, hubert2013scaling}.

Prominent classes of chemical networks exhibiting such structural properties are certain \emph{multisite phosphorylation networks} \cite{thomson2009unlimited} and the \emph{Monod--Wyman--Changeux model} of allosteric ligand binding \cite{ gunawardena2012linear, monod1965nature}.
Furthermore, symmetric structures can be found in \emph{autocatalytic recombination systems} and \emph{hypercycles}, which are extensively studied  in ``origin of life'' models to understand the origins of biological complexity \cite{craciun2022autocatalytic, eigen1971selforganization, eigen1977principle, hordijk2010autocatalytic}.

The exploitation of symmetries for finding points in symmetric semialgebraic sets has been a fruitful research direction \cite{faugere2023computing,Riener12,riener2018real}. There the symmetries of reflection groups are used to bound the number of distinct coordinates and to make computations effective. We make use of these results and connect them to steady states of symmetric CRNs.\\

In this paper, we systematically study the 
relationship between symmetric CRNs and equivariant mass-action ODEs. We prove that equivariant mass-action ODEs are always realized by symmetric CRNs, show that the graph of the dense realization is symmetric and construct a sparsest symmetric realization for rational coefficients (Section \ref{sec: 2}). Moreover, we demonstrate how Timofte's half-degree principle can be applied to reduce the search for steady states to solutions in finitely many lower dimensional subspaces (Section \ref{sec: 3}).

\section*{Notation}
For a positive integer $n$ we write $[n]:=\{1,2,\dots,n\}$, $\RR[x]:=\RR[x_1,\dots,x_n]$ and $S_n$ for the symmetric group on $[n]$. Let $\N=(\S,\C,\R)$ be a chemical reaction network with $n$ species $X_1,\dots,X_n$. We assume that $\N$ has no redundant complexes, i.e., we suppose that $A \in \mathcal{C}$ implies $A \to B \in \mathcal{R}$ or $B \to A \in \mathcal{R}$ for some $B \in \mathcal{C}$. Note that this is not a restriction, since complexes that are not present in any reaction have no influence on the mass-action dynamics of the CRN. 
In mass-action kinetics one associates to each CRN a vector function \[f(x)=(f_1(x),\dots,f_n(x))=\sum_{A \to B \in \mathcal{R}}\kappa_{A,B} x^a(b-a),\] where $a,b \in \NN_0^n$ are the \emph{coefficient vectors} of the complexes $A=\sum_{i=1}^na_iX_i$ and $B=\sum_{i=1}^nb_iX_i$, $x^a := x_1^{a_1}\cdots x_n^{a_n}$, and $\kappa_{A,B} >0$ are fixed reaction parameters. For $x^a$ we call $a$ the \emph{exponent vector}. We call $f(x)$ the \emph{mass-action ODE} of $\mathcal{N}$.
For a given CRN we will always assume that the reaction parameters are given. 
Observe that each coordinate $f_i(x)$ of $f(x)$ is a real polynomial in $x_1,\ldots,x_n$. 
A \emph{steady state} of the CRN $\mathcal{N}$ is a solution of the system of polynomial equations $f(x)=0$ in the non-negative orthant $\RR_{\geq 0}^n$. We call a steady state \emph{positive} if all of its coordinates are positive, i.e., the steady state is contained in the positive orthant $\RR_{> 0}^n$.

\section{Symmetric chemical reaction networks and \\
equivariant mass-action ODEs} \label{sec: 2}

In this section we first consider general symmetries on CRNs which are purely graph-theoretic.
Since we will rely on the notation, we reprove the well known result in equivariant network theory that the mass-action ODE of a symmetric CRN is equivariant (Proposition \ref{prop: symmetric CRN implies mass-action ODE equivariant}).
We then focus on those symmetries that arise from permutations of the species of the CRN.
For these actions we prove that any equivariant mass-action ODE is realized by a symmetric CRN (Proposition \ref{prop: mass-action ODE invariant implies symmetric CRN}). Moreover, we deduce that the sum of squares of the coordinate functions of an equivariant mass-action ODE is invariant (Corollary \ref{cor: F invariant}).
Finally, we consider symmetric realizations of given equivariant mass-action ODEs and study the densest and a sparsest realization (Proposition \ref{prop: unweighted symmetric} and Theorem \ref{thm: sparsest}).

\begin{example}[{\cite[Figure~5(b)]{craciun2022autocatalytic}}]\label{ex:dickenstein} 
The CRN $\N_1$ in Figure \ref{fig:N1} with all reaction parameters equal to $1$ has various symmetries.
    We can permute the species using the permutation $(123)$ and all complexes and reactions are preserved.
    This symmetry corresponds to the $120^\circ$-rotation.
    Additionally, it admits the graph theoretic symmetry corresponding to reflection about the horizontal axis which is no permutation on the species.
    So, on the set of species we find a $\ZZ_3$-symmetry, while the symmetry group $\operatorname{Aut}(\N_1)$ of the underlying directed graph even has $6$ elements and is isomorphic to $S_3$.
    The CRN $\N_1$ has the mass-action ODE
    \begin{equation*}
        f_{\N_1}=\begin{pmatrix}
            &-2x_1^2x_2&-x_1^2x_3&&&+2x_2^2x_3&+x_2x_3^2\\
            &&\phantom{+}x_1^2x_3&-x_1x_2^2&+2x_1x_3^2&-2x_2^2x_3&\\
            &\phantom{+}2x_1^2x_2&&+x_1x_2^2&-2x_1x_3^2&&-x_2x_3^2
        \end{pmatrix}.
    \end{equation*}
    In Figure \ref{fig:example21}, the concentrations of the species and complexes of $\N_1$ with randomly sampled reaction parameters over time are shown.
\end{example}

\begin{figure}
    \centering
    \begin{center}
        \begin{tikzpicture}
            
            \node[draw, circle, white, fill=white, minimum size=1.6mm, inner sep=0pt] (1b) at (0:1.9) {};
            \node[draw, circle, white, fill=white, minimum size=1.6mm, inner sep=0pt] (1c) at (5:1.9) {};
            \node[draw, circle, white, fill=white, minimum size=1.6mm, inner sep=0pt] (1d) at (-5:1.9) {};
            \node[draw, circle, fill=black, minimum size=1.6mm, inner sep=0pt] (1a) at (0:2) {};
            \node at (0:3) {$2X_2+X_3$};
            
            \node[draw, circle, white, fill=white, minimum size=1.6mm, inner sep=0pt] (2b) at (60:1.9) {};
            \node[draw, circle, white, fill=white, minimum size=1.6mm, inner sep=0pt] (2c) at (65:1.9) {};
            \node[draw, circle, white, fill=white, minimum size=1.6mm, inner sep=0pt] (2d) at (55:1.9) {};
            \node[draw, circle, fill=black, minimum size=1.6mm, inner sep=0pt] (2a) at (60:2) {};
            \node at (60:2.5) {$X_2+2X_3$};
            
            \node[draw, circle, white, fill=white, minimum size=1.6mm, inner sep=0pt] (3b) at (120:1.9) {};
            \node[draw, circle, white, fill=white, minimum size=1.6mm, inner sep=0pt] (3c) at (125:1.9) {};
            \node[draw, circle, white, fill=white, minimum size=1.6mm, inner sep=0pt] (3d) at (115:1.9) {};
            \node[draw, circle, fill=black, minimum size=1.6mm, inner sep=0pt] (3a) at (120:2) {};
            \node at (120:2.5) {$X_1+2X_3$};
            
            \node[draw, circle, white, fill=white, minimum size=1.6mm, inner sep=0pt] (4b) at (180:1.9) {};
            \node[draw, circle, white, fill=white, minimum size=1.6mm, inner sep=0pt] (4c) at (185:1.9) {};
            \node[draw, circle, white, fill=white, minimum size=1.6mm, inner sep=0pt] (4d) at (175:1.9) {};
            \node[draw, circle, fill=black, minimum size=1.6mm, inner sep=0pt] (4a) at (180:2) {};
            \node at (180:3) {$2X_1+X_3$};
            
            \node[draw, circle, white, fill=white, minimum size=1.6mm, inner sep=0pt] (5b) at (240:1.9) {};
            \node[draw, circle, white, fill=white, minimum size=1.6mm, inner sep=0pt] (5c) at (245:1.9) {};
            \node[draw, circle, white, fill=white, minimum size=1.6mm, inner sep=0pt] (5d) at (235:1.9) {};
            \node[draw, circle, fill=black, minimum size=1.6mm, inner sep=0pt] (5a) at (240:2) {};
            \node at (240:2.5) {$2X_1+X_2$};
            
            \node[draw, circle, white, fill=white, minimum size=1.6mm, inner sep=0pt] (6b) at (300:1.9) {};
            \node[draw, circle, white, fill=white, minimum size=1.6mm, inner sep=0pt] (6c) at (305:1.9) {};
            \node[draw, circle, white, fill=white, minimum size=1.6mm, inner sep=0pt] (6d) at (295:1.9) {};
            \node[draw, circle, fill=black, minimum size=1.6mm, inner sep=0pt] (6a) at (300:2) {};
            \node at (300:2.5) {$X_1+2X_2$};

            \node[draw, circle, white, fill=white, minimum size=1.6mm, inner sep=0pt] (7b) at (0:0.1) {};
            \node[draw, circle, white, fill=white, minimum size=1.6mm, inner sep=0pt] (7c) at (60:0.1) {};
            \node[draw, circle, white, fill=white, minimum size=1.6mm, inner sep=0pt] (7d) at (120:0.1) {};
            \node[draw, circle, white, fill=white, minimum size=1.6mm, inner sep=0pt] (7e) at (180:0.1) {};
            \node[draw, circle, white, fill=white, minimum size=1.6mm, inner sep=0pt] (7f) at (240:0.5) {};
            \node[draw, circle, white, fill=white, minimum size=1.6mm, inner sep=0pt] (7g) at (300:0.5) {};
            \node[draw, circle, fill=black, minimum size=1.6mm, inner sep=0pt] (7a) at (0:0) {};
            \node at (0,-0.3) {$X_1+X_2+X_3$};
            
            \draw[arrows = 
            {-Stealth[length=0.2cm]},style={shorten <=4pt, shorten >=4pt}] (1a) to (7a);
            \draw[arrows = 
            {-Stealth[length=0.2cm]},style={shorten <=4pt, shorten >=4pt}] (7a) to (2a);
            \draw[arrows = 
            {-Stealth[length=0.2cm]},style={shorten <=4pt, shorten >=4pt}] (3a) to (7a);
            \draw[arrows = 
            {-Stealth[length=0.2cm]},style={shorten <=4pt, shorten >=4pt}] (7a) to (4a);
            \draw[arrows = 
            {-Stealth[length=0.2cm]},style={shorten <=4pt, shorten >=15pt}] (5a) to (7a);
            \draw[arrows = 
            {-Stealth[length=0.2cm]},style={shorten <=15pt, shorten >=4pt}] (7a) to (6a);
            
            \draw[arrows = 
            {-Stealth[length=0.2cm]}] (1c.70) to (2d.350);
            \draw[arrows = 
            {-Stealth[length=0.2cm]}] (2d.250) to (1c.170);

            \draw[arrows = 
            {-Stealth[length=0.2cm]}] (2c.130) to (3d.50);
            \draw[arrows = 
            {-Stealth[length=0.2cm]}] (3d.310) to (2c.230);

            \draw[arrows = 
            {-Stealth[length=0.2cm]}] (3c.190) to (4d.110);
            \draw[arrows = 
            {-Stealth[length=0.2cm]}] (4d.10) to (3c.290);

            \draw[arrows = 
            {-Stealth[length=0.2cm]}] (4c.250) to (5d.170);
            \draw[arrows = 
            {-Stealth[length=0.2cm]}] (5d.70) to (4c.350);

            \draw[arrows = 
            {-Stealth[length=0.2cm]}] (5c.310) to (6d.230);
            \draw[arrows = 
            {-Stealth[length=0.2cm]}] (6d.130) to (5c.50);

            \draw[arrows = 
            {-Stealth[length=0.2cm]}] (1d.190) to (6c.110);
            \draw[arrows = 
            {-Stealth[length=0.2cm]}] (6c.10) to (1d.290);
        \end{tikzpicture}
    \end{center}
    \caption{The CRN $\N_1$}
    \label{fig:N1}
\end{figure}

\begin{figure}
	\begin{centering}
		\includegraphics{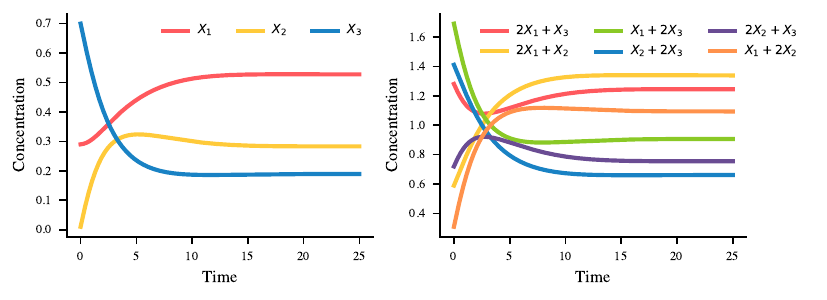}
	\end{centering}
	\caption{Simulation of the CRN $\mathcal{N}_1$ in Example 2.1 with randomly sampled reaction parameters. The left panel displays the time evolution of the monomer concentrations, while the right panel displays the evolution of the triplet concentrations. By mass conservation in the ODE system, the total concentration $x_1+x_2+x_3$ remains constant over time and is therefore omitted from the plot.}
	\label{fig:example21}
\end{figure}

We begin our analysis with the most general notion of a group action on a CRN and later we examine a more restricted setting in detail, which excludes the $S_3$-symmetry of the CRN $\N_1$ in Example \ref{ex:dickenstein}.

Recall that a \emph{left action} of a group $G$ on a set $X$ is a map $\cdot : G \times X \to X$ satisfying:
\begin{enumerate}
    \item[(i)] For all $\sigma \in G$ the map $\sigma : X \to X, x \mapsto \sigma \cdot x$ is bijective.
    \item[(ii)] For all $\sigma, \tau \in G$ we have $\sigma \cdot (\tau \cdot x) = (\sigma \tau) \cdot x$ for all $x \in X$.
    \item[(iii)] $\operatorname{id} \cdot x= x$ for all $x \in X$.
\end{enumerate}
Analogously, one can define a \emph{right action} of $G$ on $X$. Group actions are strongly connected to \emph{group representations} which are group homomorphisms $G \to \operatorname{GL}(V)$ for vector spaces $V$. In fact, any group action canonically induces a group representation on the free vector space over $X$. 

Let $\N=(\S,\C,\R)$ with $n$ species $\S=\{X_1,\dots,X_n\}$ be a CRN equipped with a left action of a group $G$ on the set of complexes $\C$ which we denote by 
\begin{align}
 \sigma \circ (a_1X_1+\ldots+a_nX_n)  \label{eq: left-action on complexes}   
\end{align} 
for $\sigma \in G$ and $a_1X_1+\ldots+a_nX_n  \in \mathcal{C}$. 
This action on complexes naturally induces an action on the set of coefficient vectors of complexes, as a subset of $\NN_0^n$, in the following way: If $\sigma \circ A= B$ for complexes $A, B \in \C$ with coefficient vectors $a, b \in \NN_0^n$, we write $\sigma \circ a=b$. The action \eqref{eq: left-action on complexes} extends naturally to an action on the set $\widetilde{\R}=\setcond{A\to B}{A,B\in\C}$ of all possible permuted reactions of $\N$ by
\begin{align}
\sigma \circ \left( A \to B \right) = \left( \sigma \circ A\right) \to  \left( \sigma \circ B \right) \label{eq: left-action on reactions}    
\end{align} 
for $\sigma \in G$ and complexes $A, B \in \C$. The set $\widetilde{\R}$ is in general a superset of $ \R$, since we only assumed that $\C$ is closed under the action of $G$. 
The purpose of this left action of $G$ is to reflect the symmetries of the directed graph of the CRN encoded by the group $G$. Moreover, $\sigma (\mathcal{N})$ is the CRN with species $\mathcal{S}$, complexes $\mathcal{C}$ and reactions $\sigma \circ (A \to B)$ for all $A \to B \in \mathcal{R}$. 

We want to extend this action of $G$ on $\C$ also to an action on the associated mass-action ODE. 
We denote the set of all \emph{polynomials supported on $\C$} by 
\[\RR[x]_\C := \left\{\sum_{A \in \C} c_A x^a \mid c_A \in \RR\right\},\]
where $a$ is the coefficient vector of the complex $A$. In particular, $\RR[x]_\C $ contains all real polynomials that are linear combinations of monomials having exponents in $\C$.
We define a right action of $G$ on $\RR[x]_\C$ by
\begin{align}\label{eq: right action}
    \sum_{A \in \C} c_A x^a  \bullet \sigma = \sum_{A \in \C} c_A x^{\sigma^{-1}(a)} =   \sum_{A \in \C} c_{\sigma(A)} x^{a}
\end{align} 
for $\sigma\in G$.
Although the definition of this right action may seem artificial, Proposition \ref{prop: symmetric CRN implies mass-action ODE equivariant} illustrates that it is in fact the canonical choice.

Moreover, for $\sigma \in G$
and the mass-action ODE $f(x)$ of a 
CRN $\N=(\S,\C,\R)$, we define
\begin{align} \label{eq:second left action}
    \sigma \ast \sum_{A \to B \in \R} \kappa_{A,B} x^a (b-a) = \sum_{A \to B \in \R} \kappa_{A,B}x^a(\sigma (b)-\sigma(a)).
\end{align}
This depends on the CRN $\N$ and not only on $f(x)$, since a given mass-action ODE could be realized by several CRNs.
\begin{definition}\label{def:symmetric and equivariant crn}
A CRN $\N=(\S,\C,\R)$ with mass-action ODE $f_{\N}(x)$ is \emph{symmetric} with respect to the group $G$ acting on $\C$ if the set of reactions is closed under the action of $G$, i.e., $\left(\sigma \circ A\right) \to \left(\sigma \circ B\right) \in \R $ for all $A \to B \in \R$ and $\sigma \in G$, and the reaction parameters $\kappa_{A,B}$ are constant on orbits of reactions. \\
We call the mass-action ODE $f_\N(x)$ \emph{equivariant} if $f_\N(x) \bullet \sigma = \sigma \ast f_\N(x)$ for all $\sigma \in G$.
\end{definition}

We illustrate the different actions $\circ, \bullet$ and $\ast$ being used in Definition \ref{def:symmetric and equivariant crn} in the case of the symmetric group.
\begin{example}\label{ex:symmetric group action}
    Let $S_n$ act on the set of species $\S = \{X_1,\ldots,X_n\}$ of a CRN $(\S,\C,\R)$ with mass-action ODE $f=(f_1,\dots,f_n)$ by permutation. This naturally induces an action on $\C$ and $\R$. The action $\circ$ on complexes is given by 
    \[ \sigma \circ (a_1X_1+\ldots +a_nX_n) = a_1X_{\sigma(1)}+\ldots+a_nX_{\sigma(n)}=a_{\sigma^{-1}(1)}X_1+\ldots+a_{\sigma^{-1}(n)}X_n.\]
    Moreover, we consider the actions $\bullet$ and $\ast$ on $\RR[x]^n$.
    For $\sigma \in S_n$ we have
    \begin{align*}
        f(x) \bullet \sigma  & = \sum_{A \to B \in \R}\kappa_{A,B} x_{1}^{a_{\sigma (1)}}\cdots x_{n}^{a_{\sigma(n)}}(b-a) = \sum_{A \to B \in \R}\kappa_{A,B} x_{\sigma^{-1}{(1)}}^{a_1}\cdots x_{\sigma^{-1}(n)}^{a_n}(b-a) \\
        & =f(x_{\sigma^{-1}(1)},\dots,x_{\sigma^{-1}(n)})
    \end{align*}
    and
    \[\sigma \ast f(x) = \sum_{A \to B \in \R} \kappa_{A,B} x^a (\sigma(b)-\sigma(a)) = (f_{\sigma^{-1}(1)}(x),\dots,f_{\sigma^{-1}(n)}(x)). \]
    Thus, $\bullet$ permutes the variables of $f$, while $\ast$ permutes the coordinate functions of $f$. 
\end{example}

Recall that the CRN $\N_1$ in Example \ref{ex:dickenstein} has the symmetry of the cyclic group $\ZZ_3$ acting on the set of species and the symmetry of the symmetric group $S_3$ acting on its directed graph.
\begin{example}\label{ex: symmetries of ODE} 
We continue with Example \ref{ex:dickenstein}.
Trivially, $\N_1$ is $\operatorname{Aut}(\N_1)$-symmetric.
We find that $f_{\N_1}$ is $\operatorname{Aut}(\N_1)$-equivariant.
To see this, it is enough to check the equivariance on the generators $(123)$, corresponding to the permutation of the species, and $s$, corresponding to the reflection symmetry.
Indeed, we get that
\begin{equation*}
    (123)\ast f_{\N_1}=\begin{pmatrix}
        &\phantom{+}2x_1^2x_2&&+x_1x_2^2&-2x_1x_3^2&&-x_2x_3^2\\
        &-2x_1^2x_2&-x_1^2x_3&&&+2x_2^2x_3&+x_2x_3^2\\
        &&\phantom{+}x_1^2x_3&-x_1x_2^2&+2x_1x_3^2&-2x_2^2x_3&
    \end{pmatrix}=f_{\N_1}\bullet(123),
\end{equation*}
and
\begin{equation*}
    s\ast f_{\N_1}=\begin{pmatrix}
        &&-x_1^2x_3&+x_1x_2^2&-2x_1x_3^2&+2x_2^2x_3&\\
        &\phantom{+}2x_1^2x_2&+x_1^2x_3&&&-2x_2^2x_3&-x_2x_3^2\\
        &-2x_1^2x_2&&-x_1x_2^2&+2x_1x_3^2&&+x_2x_3^2
    \end{pmatrix}=f_{\N_1}\bullet s.
\end{equation*}

\end{example}
We recover the following result from equivariant dynamical systems \cite{golubitsky2003symmetry} which we derive here directly from Definition \ref{def:symmetric and equivariant crn}.
\begin{proposition}\label{prop: symmetric CRN implies mass-action ODE equivariant}
If $\mathcal{M}$ is a $G$-symmetric CRN with mass-action ODE $f(x)$ then $f(x)$ is $G$-equivariant. 
\end{proposition}
\begin{proof}
 For $\sigma \in G$ we obtain
    \begin{align*}
         f(x) \bullet \sigma      & = \sum_{A\to B\in \mathcal{R}}\kappa_{A,B} x^{\sigma^{-1}(a)}(b-a)\\
        &=\sum_{A\to B\in \mathcal{R}}\kappa_{\sigma \circ A,\sigma \circ B}x^{\sigma^{-1}(\sigma(a))}(\sigma(b)-\sigma(a))\\
        &=\sum_{A\to B\in \mathcal{R}}\kappa_{A,B}x^{a}\left(\sigma(b)-\sigma(a)\right)\\
        &=\sigma\ast f(x),
    \end{align*}
    where the second equality holds because $\sigma \in G$ acts bijectively on $\mathcal{R}$ and the reaction parameters are constant on every orbit by assumption.
\end{proof}

In the following, we address the question whether every equivariant vector function in $\RR[x]^n$, which is the mass-action ODE of a CRN, can be realized by a symmetric CRN. 
In fact, there can be infinitely many CRNs with the same mass-action ODE (see e.g. \cite[Section~4]{craciun2008identifiability}). 

\begin{example}\label{ex: asymmetrix CRN}
    The mass-action ODE $f_{\N_1}$ is $\operatorname{Aut}(\N_1)$-equivariant.
    However, not every CRN with mass-action ODE $f_{\N_1}$ is symmetric with respect to $\operatorname{Aut}(\N_1)$, or with respect to any non-trivial group.
    Consider for example the CRN $\N_2$ with all reaction parameters equal to $1$:
    \begin{center}
    \begin{tikzpicture}[dot/.style={circle, fill, inner sep=1.8pt},inner/.style={shorten <=5pt, shorten >=6pt}]
    
    \node[dot, label={0:{$2X_2+X_3$}}]     (1) at (0:2.2)     {};
    \node[dot, label={60:{$X_2+2X_3$}}]    (2) at (60:2.2)    {};
    \node[dot, label={120:{$X_1+2X_3$}}]   (3) at (120:2.2)   {};
    \node[dot, label={120:{$2X_1+X_3$}}]   (4) at (180:2.2)   {};
    \node[dot, label={240:{$2X_1+X_2$}}]   (5) at (240:2.2)   {};
    \node[dot, label={300:{$X_1+2X_2$}}]   (6) at (300:2.2)   {};
    \node[dot, label={270:{$X_1+X_2+X_3$}}]  (7) at (0,0)     {};
        
    \draw[-{Stealth[length=0.3cm]}] (1) -- (6);  
    \draw[-{Stealth[length=0.3cm]}] (3) -- (4);   
    \draw[-{Stealth[length=0.3cm]}] (5) -- (6);  
    
    \draw[-{Stealth[length=0.3cm]}] (2) -- (7);   
    \draw[-{Stealth[length=0.3cm]}] (4) -- (7);   
    \draw[-{Stealth[length=0.3cm]}] (6) -- (300:0.7);   
    
    \draw[-{Stealth[length=0.3cm]},inner] (1) edge[bend right=22] (3);   
    \draw[-{Stealth[length=0.3cm]},inner] (3) edge[bend left=12] (1); 
    \draw[-{Stealth[length=0.3cm]},inner] (5) edge[bend left=22] (3);

    \end{tikzpicture}
    \end{center}

    We find $f_{\N_1}=f_{\N_2}$, where $f_{\N_2}$ denotes the mass-action ODE of $\N_2$.
    
   \end{example}

We defined with $\bullet $ in \eqref{eq: right action} and $\ast$ in \eqref{eq:second left action} two operations of a group $G$ on mass-action ODEs which strongly depend on the action of $G$ on the CRN. 
We want to use a similar operation for general vector functions. 
For this, we need a well-defined action on the integer coefficient vectors and on the monomials, but also on any complex of an associated CRN. 
Therefore we restrict in the following to linear group actions which map $\NN_0^n$ onto $\NN_0^n$. 
The following lemma shows that this actually means that we are restricting to subgroups of the permutation group of the set of species.

\begin{lemma} \label{lem:symmetric group is natural}
Any invertible matrix $A \in \RR^{n \times n}$ which maps $\NN_0^n$ onto $\NN_0^n$ is a permutation matrix.
\end{lemma}
\begin{proof}
Let $\mathrm{e}_1,\ldots,\mathrm{e}_n$ denote the standard basis vectors in $\RR^n$. Then $A\mathrm{e}_i = A_{\bullet i}  \in \operatorname{Im}(A)$ must be a vector in $\NN_0^n$. 
Moreover, for every $1 \leq j \leq n$ there must exist $\lambda_1,\ldots,\lambda_n \in \NN_0$ such that $\mathrm{e}_j = A(\lambda_1\mathrm{e}_1+\ldots+\lambda_n \mathrm{e}_n)=\sum_{i=1}^n \lambda_i A_{\bullet i}$ which proves that $A$ is indeed a permutation of the columns of the identity matrix $\begin{pmatrix}
    \mathrm{e}_1 | \mathrm{e}_2 | \cdots | \mathrm{e}_n
\end{pmatrix}$.
\end{proof}
In the remaining part of the paper we restrict to permutation groups which act as in Example \ref{ex:symmetric group action}.
\begin{definition}\label{def:2}
For $G\subset S_n$ we call a vector function $h(x) \in \RR[x]^n$ $G$-\emph{equivariant}, if $h(x) \bullet \sigma = \sigma \ast h(x)$ for all $\sigma \in G$. 
\end{definition}
Observe that if $h(x)$ is the mass-action ODE of a CRN then Definition \ref{def:2} coincides with the second part of Definition \ref{def:symmetric and equivariant crn}.

The following result will be the key to Section \ref{sec: 3}.

\begin{corollary}\label{cor: F invariant}
    Let $G \subset S_n$ and $h=(h_1,\dots,h_n)\in \RR[x]^n$ be $G$-equivariant. Then $H=h_1^2+\dots+h_n^2 \in \RR[x]$ is $G$-invariant, i.e., $H(x_{\sigma(1)},\dots,x_{\sigma(n)})=H(x)$ for all $\sigma \in G$.
\end{corollary}
\begin{proof}
    Since $h$ is $G$-equivariant, we have $h(x) \bullet \sigma^{-1} = \sigma^{-1} \ast h(x)$ for all $\sigma \in G$.
    As illustrated in Example \ref{ex:symmetric group action}, $h(x) \bullet \sigma^{-1}$ corresponds to evaluating each function $h_i$ in the vector function $h$ at $(x_{\sigma(1)},\dots, x_{\sigma(n)})$ and $\sigma^{-1} \ast h(x) = (h_{\sigma(1)}(x),\dots,h_{\sigma(n)}(x))$. Therefore, evaluating $H$ at these coordinates yields
    \[ H(x_{\sigma(1)},\dots,x_{\sigma(n)}) = \sum_{i=1}^n h_i(x_{\sigma(1)},\dots,x_{\sigma(n)})^2 = \sum_{i=1}^n h_{\sigma(i)}(x)^2 = H(x).\qedhere\]  
\end{proof}

The converse of Proposition \ref{prop: symmetric CRN implies mass-action ODE equivariant} is in general not true (see Example \ref{ex: asymmetrix CRN}).
Which polynomial vector functions are realizable as the mass-action ODE of a CRN is characterized by the \emph{Hungarian lemma} \cite{hars1981inverse}.
Also, if we only consider the realizable vector functions, then the converse remains false as is shown by Example \ref{ex: asymmetrix CRN}.
It is however true that a symmetric realization always exists.
We show this in Proposition \ref{prop: mass-action ODE invariant implies symmetric CRN}.

\begin{lemma}\label{lem: actonCRN}
    Let $G$ be a group acting on the mass-action ODE $f_\mathcal{M}$ of some CRN $\mathcal{M}$.
    Then for all $\sigma\in G$, we have
    \begin{equation*}
        f_{\sigma(\mathcal{M})}=\sigma\ast f_\mathcal{M}\bullet\sigma^{-1}.
    \end{equation*}
    In particular, if $f_{\mathcal{M}}$ is $G$-equivariant, then $f_{\sigma(\mathcal{M})}=f_\mathcal{M}$.
\end{lemma}
\begin{proof}
    For every reaction $A\to B\in\mathcal{R}_\mathcal{M}$ with reaction parameter $\kappa_{A,B}$, the CRN $\sigma(\mathcal{M})$ contains the reaction $\sigma(A)\to\sigma(B)$ with the same reaction parameter $\kappa_{A,B}$.
    Consequently, the mass-action ODE of the new CRN $\sigma(\mathcal{M})$ reads
    \begin{align*}
        f_{\sigma(\mathcal{M})}&=\sum_{A\to B\in\mathcal{R}_\mathcal{M}}\kappa_{A,B}x^{\sigma(a)}(\sigma(b)-\sigma(a))\\
        &=\left(\sum_{A\to B\in\mathcal{R}_\mathcal{M}}\kappa_{A,B}x^a(\sigma(b)-\sigma(a))\right)\bullet\sigma^{-1}\\
        &=\sigma\ast\left(\sum_{A\to B\in\mathcal{R}_\mathcal{M}}\kappa_{A,B}x^a(b-a)\right)\bullet\sigma^{-1}=\sigma\ast f_{\mathcal{M}}\bullet\sigma^{-1}.
    \end{align*}
    If $f_\mathcal{M}$ is $G$-equivariant, we can simplify this:
    \begin{equation*}
        f_{\sigma(\mathcal{M})}=\sigma\ast f_\mathcal{M}\bullet\sigma^{-1}=f_\mathcal{M}\bullet\sigma\bullet\sigma^{-1}=f_\mathcal{M}.\qedhere
    \end{equation*}
\end{proof}

We call the graph $(\mathcal{C},\mathcal{R})$ of a CRN its \emph{unweighted graph}, i.e., the reaction parameters are ignored. 
In \cite[Theorem 3.1 and Theorem 3.2]{szederkenyi2010maximal}, the authors proved that, given a realizable vector function and a set of possible complexes, the unweighted graph of any of its realizations is a subgraph of the unweighted graph of the unique \emph{dense} realization.
Consequently, the dense realization is the realization of the given mass-action ODE with the maximal number of reactions with respect to the fixed set of possible complexes. 

\begin{proposition}\label{prop: unweighted symmetric}
    Let $G\subset S_n$ and let $f(x)\in\RR[x]^n$ be a realizable $G$-equivariant vector function. 
    Then the unweighted graph of the dense realization, with respect to a fixed $G$-invariant set of complexes, is $G$-symmetric.
\end{proposition}
\begin{proof}
    Let $\mathcal{M}$ be any realization.
    Then Lemma \ref{lem: actonCRN} shows that $\sigma (\mathcal{M})$ is a realization as well.
    From \cite[Theorem 3.1]{szederkenyi2010maximal}, it follows that the unweighted graph of the dense realization $\mathcal{D}$ contains in particular the unweighted graph of the realization $\sigma(\mathcal{D})$ for all $\sigma\in G$ and thus that $\mathcal{D}$ is $G$-symmetric modulo the reaction parameters.
\end{proof}

In this form, Proposition \ref{prop: unweighted symmetric} is not constructive and only gives us symmetry of the unweighted graph of the dense realization.
However, symmetry of a CRN is also a condition on the reaction parameters.
The following is a constructive result.

\begin{proposition}\label{prop: mass-action ODE invariant implies symmetric CRN}
Let $G \subset S_n$.  
If $h(x)\in \RR[x]^n$ is $G$-equivariant and there exists a CRN inducing $h(x)$ as mass-action ODE, then there exists a $G$-symmetric CRN inducing $h(x)$.  
\end{proposition}
\begin{proof}
We aim to construct a $G$-symmetric CRN $\mathcal{N}$ whose mass-action ODE is $h(x)$.
Let $\mathcal{M}=(\S_\M,\C_\M,\R_\M)$ be any CRN with mass-action ODE $h(x)$.
For each $\sigma \in G$ and reaction $A\to B \in \R_\M$ with reaction parameter $\kappa^\M_{A,B}$, we generate a reaction $\sigma \circ (A\to B)$ in $\N$ with reaction parameter $\frac{\kappa^\mathcal{M}_{A,B}}{|G|}$. If multiple reactions in $\mathcal{M}$ and permutations in $G$ generate the exact same reaction in $\mathcal{N}$, their reaction parameters are summed.

Let $g(x)$ denote the mass-action ODE of $\mathcal{N}$. Summing over all $\sigma \in G$ and all reactions in $\mathcal{M}$, we find:
\begin{align*}
    g(x) &= \sum_{\sigma \in G} \sum_{A \to B \in \mathcal{R}_\mathcal{M}} \frac{1}{|G|} \kappa^\mathcal{M}_{A,B} x^{\sigma(a)} (\sigma(b) - \sigma(a)) \\
    &= \frac{1}{|G|} \sum_{\sigma \in G} \sigma \ast \left( \sum_{A \to B \in \mathcal{R}_\mathcal{M}} \kappa^\mathcal{M}_{A,B} x^a (b-a) \right) \bullet \sigma^{-1} \\
    &= \frac{1}{|G|} \sum_{\sigma \in G} \sigma \ast h(x) \bullet \sigma^{-1}
\end{align*}
Because $h(x)$ is $G$-equivariant, we have $h(x) \bullet \sigma = \sigma \ast h(x)$, which implies $\sigma \ast h(x) \bullet \sigma^{-1} = h(x)$. Therefore:
\[g(x) = \frac{1}{\vert{}G\vert{}} \sum_{\sigma \in G} h(x) = h(x).\]
Thus, the symmetric CRN $\mathcal{N}$ realizes the mass-action ODE $h(x)$.
\end{proof}
Note that since reaction parameters are always positive, in this construction they can never cancel out.
Consequently, the symmetric CRN that is constructed contains the original CRN.
Thus, starting with a dense realization and applying the construction from Proposition \ref{prop: mass-action ODE invariant implies symmetric CRN} will give us a realization on the exact same graph, but possibly with different reaction parameters.
In particular, for every realizable $G$-equivariant vector function, there is a symmetric realization on the dense CRN graph.
The following example shows that Proposition \ref{prop: mass-action ODE invariant implies symmetric CRN} does not always yield a dense realization.

\begin{example}
    Consider the mass-action ODE 
    \begin{equation*}
        f_{\N_1}=\begin{pmatrix}
            &-2x_1^2x_2&-x_1^2x_3&&&+2x_2^2x_3&+x_2x_3^2\\
            &&\phantom{+}x_1^2x_3&-x_1x_2^2&+2x_1x_3^2&-2x_2^2x_3&\\
            &\phantom{+}2x_1^2x_2&&+x_1x_2^2&-2x_1x_3^2&&-x_2x_3^2
        \end{pmatrix}
    \end{equation*}
    from Example \ref{ex:dickenstein}.
    It has various symmetries as we have checked in Example \ref{ex: symmetries of ODE}.
    Since we start with an equivariant mass-action ODE we can also allow for a group action that does not permute the species in contrast to the setting of Proposition \ref{prop: mass-action ODE invariant implies symmetric CRN}. 
    So we work with $\operatorname{Aut}(\N_1)\cong S_3$.
    In Example \ref{ex: asymmetrix CRN} we have seen an asymmetric realization of $f_{\N_1}$.
    We apply the constructive method used in the proof of Proposition \ref{prop: mass-action ODE invariant implies symmetric CRN} to construct a CRN $\N_3$ that induces $f_{\N_1}$.
    The method yields:
    \begin{center}
    \begin{tikzpicture}[dot/.style={circle, fill, inner sep=1.8pt},inner/.style={shorten <=5pt, shorten >=6pt}]
    
    \node[dot, label={0:{$2X_2+X_3$}}]     (1) at (0:2.7)     {};
    \node[dot, label={60:{$X_2+2X_3$}}]    (2) at (60:2.7)    {};
    \node[dot, label={120:{$X_1+2X_3$}}]   (3) at (120:2.7)   {};
    \node[dot, label={180:{$2X_1+X_3$}}]   (4) at (180:2.7)   {};
    \node[dot, label={240:{$2X_1+X_2$}}]   (5) at (240:2.7)   {};
    \node[dot, label={300:{$X_1+2X_2$}}]   (6) at (300:2.7)   {};
    \node[dot, label={270:{$X_1+X_2+X_3$}}]  (7) at (0,0)     {};
    
    \draw[-{Stealth[length=0.3cm]},dashed] (1) -- (2);   
    \draw[-{Stealth[length=0.3cm]},dashed] (3) -- (2);   
    \draw[-{Stealth[length=0.3cm]},dashed] (3) -- (4);   
    \draw[-{Stealth[length=0.3cm]},dashed] (5) -- (4);   
    \draw[-{Stealth[length=0.3cm]},dashed] (5) -- (6);  
    \draw[-{Stealth[length=0.3cm]},dashed] (1) -- (6);   
    
    \draw[-{Stealth[length=0.3cm]}] (2) -- (7);   
    \draw[-{Stealth[length=0.3cm]}] (4) -- (7);   
    \draw[-{Stealth[length=0.3cm]}] (6) -- (300:0.7);   
    
    \draw[-{Stealth[length=0.3cm]},dashed,inner] (1) edge[bend right=12] (3);  
    \draw[-{Stealth[length=0.3cm]},dashed,inner] (3) edge[bend left=22]  (1);  
    
    \draw[-{Stealth[length=0.3cm]},dashed,inner] (3) edge[bend right=12] (5);  
    \draw[-{Stealth[length=0.3cm]},dashed,inner] (5) edge[bend left=22]  (3);  
    
    \draw[-{Stealth[length=0.3cm]},dashed,inner] (5) edge[bend right=12] (1);  
    \draw[-{Stealth[length=0.3cm]},dashed,inner] (1) edge[bend left=22]  (5);  
    
    \end{tikzpicture}
    \end{center}
    Here, the dashed arrows indicate a reaction parameter $1/2$, whereas the solid arrows towards the center have reaction parameter $1$.
    Note that $\N_3$ is not a dense realization of $f_{\N_1}$.
    Thus, the symmetric realization obtained from the method of Proposition \ref{prop: mass-action ODE invariant implies symmetric CRN} depends on the choice of CRN we start with.
    There are more symmetric realizations, in fact, $\N_3$ is, in the unweighted sense, the union of two $\ZZ_3$-symmetric realizations of $f_{\N_1}$.
    
    \begin{minipage}{0.5\textwidth}
    \begin{center}
    \begin{tikzpicture}[dot/.style={circle, fill, inner sep=1.8pt},inner/.style={shorten <=5pt, shorten >=6pt}]
    
    \node[dot, label={300:{$2X_2+X_3$}}]     (1) at (0:2.2)     {};
    \node[dot, label={60:{$X_2+2X_3$}}]    (2) at (60:2.2)    {};
    \node[dot, label={120:{$X_1+2X_3$}}]   (3) at (120:2.2)   {};
    \node[dot, label={180:{$2X_1+X_3$}}]   (4) at (180:2.2)   {};
    \node[dot, label={240:{$2X_1+X_2$}}]   (5) at (240:2.2)   {};
    \node[dot, label={300:{$X_1+2X_2$}}]   (6) at (300:2.2)   {};
    \node[dot, label={270:{$X_1+X_2+X_3$}}]  (7) at (0,0)     {};
    
    \draw[-{Stealth[length=0.3cm]}] (1) -- (2);   
    \draw[-{Stealth[length=0.3cm]}] (3) -- (4);   
    \draw[-{Stealth[length=0.3cm]}] (5) -- (6);  
    
    \draw[-{Stealth[length=0.3cm]}] (2) -- (7);   
    \draw[-{Stealth[length=0.3cm]}] (4) -- (7);   
    \draw[-{Stealth[length=0.3cm]}] (6) -- (300:0.7);   
    
    \draw[-{Stealth[length=0.3cm]},inner] (3) edge[bend left=22]  (1);  
    
    \draw[-{Stealth[length=0.3cm]},inner] (5) edge[bend left=22]  (3);  
     
    \draw[-{Stealth[length=0.3cm]},inner] (1) edge[bend left=22]  (5);  
    
    \end{tikzpicture}
    \end{center}
    \end{minipage}
    \begin{minipage}{0.5\textwidth}
    \begin{center}
    \begin{tikzpicture}[dot/.style={circle, fill, inner sep=1.8pt},inner/.style={shorten <=5pt, shorten >=6pt}]
    
    \node[dot, label={0:{$2X_2+X_3$}}]     (1) at (0:2.2)     {};
    \node[dot, label={60:{$X_2+2X_3$}}]    (2) at (60:2.2)    {};
    \node[dot, label={120:{$X_1+2X_3$}}]   (3) at (120:2.2)   {};
    \node[dot, label={120:{$2X_1+X_3$}}]   (4) at (180:2.2)   {};
    \node[dot, label={240:{$2X_1+X_2$}}]   (5) at (240:2.2)   {};
    \node[dot, label={300:{$X_1+2X_2$}}]   (6) at (300:2.2)   {};
    \node[dot, label={270:{$X_1+X_2+X_3$}}]  (7) at (0,0)     {};
    
    \draw[-{Stealth[length=0.3cm]}] (3) -- (2);     
    \draw[-{Stealth[length=0.3cm]}] (5) -- (4);   
    \draw[-{Stealth[length=0.3cm]}] (1) -- (6);   
    
    \draw[-{Stealth[length=0.3cm]}] (2) -- (7);   
    \draw[-{Stealth[length=0.3cm]}] (4) -- (7);   
    \draw[-{Stealth[length=0.3cm]}] (6) -- (300:0.7);   
    
    \draw[-{Stealth[length=0.3cm]},inner] (1) edge[bend right=22] (3);   
    
    \draw[-{Stealth[length=0.3cm]},inner] (3) edge[bend right=22] (5);   
    
    \draw[-{Stealth[length=0.3cm]},inner] (5) edge[bend right=22] (1);  
    
    \end{tikzpicture}
    \end{center}
    \end{minipage}
    Here, in both versions, the reaction parameters are $1$ for every reaction.
    In the union, we use half of every reaction parameter, and add the reaction parameters for reactions that appear in both CRNs.
\end{example}

We have seen that symmetric realizations of equivariant mass-action ODEs are not unique, and that the unique dense realization is a symmetric graph.
The remaining part of Section \ref{sec: 2} will be dedicated to constructing a sparsest symmetric realization, that is, a symmetric realization with the minimum number of reactions.
For this, let $\M_f=(\S_f,\C_f,\R_f)$ with $\S_f=\{X_1,\ldots,X_n\}$ be a sparsest $G$-symmetric realization of the $G$-equivariant mass-action ODE $f$ for some $G\subset S_n$.
Note first that the coefficient vector of every reactant of reactions in $\R_f$ must occur as exponent in $f$.
In fact, the contributions to $f$ from reactions starting at any other complex must cancel out, so deleting them from the CRN does not change the mass-action ODE, but reduces the number of reactions in the CRN.
The following lemma proves that $\M_f$ can be chosen such that whenever a symmetry in $G$ fixes a reactant, it also fixes the corresponding product.

\begin{lemma}\label{lem:at most two reactions}
    Let $\M$ be a $G$-symmetric CRN and let $A\xlongrightarrow{\kappa_{A,B}}B$ be a reaction in $\M$ where the coefficient vector of $A$ is $a$ and the stabilizer of $a$ in $G$ is $G_a$.
    Then the $|G_aB|$ reactions $A\xlongrightarrow{\kappa_{A,B}}\sigma (B)$ with $\sigma \in G_a$ in $\M$ can be replaced by at most two reactions, both of which are fixed by $G_a$, such that the mass-action ODE stays unchanged.
\end{lemma}
\begin{proof}
    Let $b$ be the coefficient vector of $B$.
    Then $G_ab$ is the set of coefficient vectors of complexes in $G_aB$.
    If $|G_aB|=1$, nothing is to be shown.
    For $|G_aB|>1,$ set 
    $$
        d=\sum_{b'\in G_ab}(b'-a)\in\ZZ^n.
    $$
    Denote by $J_-$ the coordinates corresponding to negative entries of $d$, i.e., $J_-=\setcond{i}{d_i<0}$, and set $L=\max(\setn{0}\cup\setcond{-d_i}{i\in J_-})$ and $u=\mathbf{1}_{J_-}$.
    Lastly, define 
    $$
        b^-=a-u,\qquad b^+=a+d+Lu.
    $$
    Note that both $b^-$ and $b^+$ are non-negative: $b^+$ is the sum of the two non-negative vectors $a$ and $d+Lu$, whereas the only potentially negative entries in $b^-$ are those corresponding to the coordinates in $J_-$.
    However, for $i\in J_-$, we have that $a_i\geq 1$.
    Otherwise, $b'_i-a_i\geq0$ for all $b'\in G_ab$, and thus $d_i\geq0$ contradicting the definition of $J_-$.\\

    Note further that both $a$ and $d$, and thus also $u$ are fixed by $G_a$.
    Consequently, also $b^-$ and $b^+$ are fixed by $G_a$.
    Denote by $B^-$ and $B^+$ the complexes with coefficient vectors $b^-$ and $b^+$, respectively.
    Since $d=(b^+-a)+L(b^--a)$, we get that the two reactions
    $$
        A\xlongrightarrow{\kappa_{A,B}\cdot L}B^-,\qquad A\xlongrightarrow{\kappa_{A,B}}B^+,
    $$
    together have the same contributions to the mass-action ODE as $G_a(A\xlongrightarrow{\kappa_{A,B}}B)$. 
    Note that if $d$ has no negative entries, then $b^-=a$ and $\kappa_{A,B}\cdot L=0$, and we can omit the reaction $A\xlongrightarrow{0} A$.
\end{proof}

Let $A\in\C_f$ be a reactant of a reaction in $\R_f$ with coefficient vector $a\in\NN_0^n$.
Denote the coefficient of $x^a$ in $f$ by $c_a\in\RR^n$.
Next, denote the stabilizer of $a$ in $G$ by $G_a$, as in Lemma \ref{lem:at most two reactions}.
Then the action of $G_a$ on $[n]$ gives a partition $\S_f=J_1\sqcup J_2\sqcup\ldots\sqcup J_m$.
As the action of $G_a$ fixes both $a$ and $c_a$, their coordinates corresponding to one set $J_i$ are constant on each $J_i$.
In the following, we may therefore only record one entry for each $J_i$ and by that identify vectors fixed by $G_a$, such as $a$ and $c_a$, with vectors in $\ZZ^m$ or $\RR^m$. 
In the following proposition, we will count the number of reactions in $\M_f$.
For this, $k(a,c_a)$ denotes the minimum number of reactions necessary to obtain the monomial $c_a\cdot x^a$ in $f$, i.e.,
$$
    k(a,c_a)=\min\setcond{k}{c_a=\sum_{i=1}^k \kappa_i(b_i-a), \quad \kappa_i>0, \quad b_i\in\NN_0^m\setminus\setn{a}}.
$$

\begin{proposition}\label{prop:sparsest}
    Let $f\in\RR[x]^n$ be realizable and $G$-equivariant and denote by $r_G(f)$ the minimum number of reactions in a $G$-symmetric CRN with mass-action ODE $f$.
    Then,
    $$
        r_G(f)=\sum_{a\in S/G}|G\cdot a|k(a,c_a),
    $$
    where $S$ denotes the set of vectors that occur as exponent in $f$ and $c_a$ denotes the coefficient of $x^a$ in $f$.
\end{proposition}
\begin{proof}
We count the number of reactions in our sparsest symmetric realization $\M_f$ of $f$.
Different vectors in $S/G$ do not interact, so we can count them independently.
Given a reaction $A\to B$ in $\R_f$, by Lemma \ref{lem:at most two reactions}, the entire reaction is fixed by $G_a$. 
Consequently, the $G$-orbit of the reaction consists of $|G\cdot a|=\frac{|G|}{|G_a|}$ reactions.    
As $\R_f$ is minimal, we have exactly $k(a,c_a)$ reactions with source $A$, so every source in the orbit of $A$ appears $k(a,c_a)$ times.
\end{proof}

We are now prepared to construct $\M_f$ explicitly. 
However, we only do so for mass-action ODEs $f$ with rational coefficients.
For an exponent vector $a$ and its rational coefficient $c_a\in\QQ^m$ there is a unique primitive $p_a\in\ZZ^m$ and a unique $\lambda_a\in\QQ_{>0}$ such that $c_a=\lambda_ap_a$.

\begin{theorem}\label{thm: sparsest}
    Let $f\in\QQ[x]^n$ be realizable and $G$-equivariant.
    For every representative $a\in S/G$, set
    $$
        k_a=\begin{cases}
            1,\quad a+p_a\in\NN_0^m,\\
            2,\quad a+p_a\not\in\NN_0^m.
        \end{cases}
    $$
    Then $k(a,c_a)=k_a$ and a $G$-symmetric CRN with mass-action ODE $f$ and $r_G(f)$ reactions can be constructed explicitly.
\end{theorem}
\begin{proof}
    Suppose first that $k(a,c_a)=1$, that is, $c_a=\lambda_ap_a=\kappa(b-a)$ for some $\kappa>0$ and some $b\in\NN_0^m\setminus\setn{a}$.
    Since $p_a$ is primitive and $(b-a)$ is an integer vector, we get that $(b-a)=tp_a$ for an integer $t$, and as both $\kappa>0$ and $\lambda_a>0$, we also have $t>0$.
    Thus, $tp_a+a=b\geq0$.
    Now if, for some $i$, $p_{a,i}\geq0$, then $a_i+p_{a,i}\geq0$.
    And if $p_{a,i}<0$, then $a_i+p_{a,i}\geq a_i+tp_{a,i}\geq0$.
    Consequently $a+p_a\geq0$ and thus $k(a,c_a)=k_a$.\\

    Conversely, assume that $a+p_a\geq0$.
    Denote by $A$ and $P$ the complexes with coefficient vectors $a$ and $a+p_a$, respectively.
    Then the reaction $A\xlongrightarrow{\lambda_a}P$ realizes the monomial $c_ax^a$, so again $k(a,c_a)=k_a$.\\

    It is left to show that in every other case, that is, when $a+p_a$ has at least one negative entry, two reactions suffice to realize $c_ax^a$.
    For this, similarly to Lemma \ref{lem:at most two reactions}, let $J_-=\setcond{i}{p_{a,i}<0}$, $L=\max\setcond{-p_{a,i}}{i\in J_-}$ and $u=\mathbf{1}_{J_-}$.
    As in Lemma \ref{lem:at most two reactions}, the two vectors $b^-=a-u$ and $b^+=a+p_a+Lu$ are non-negative, so let them be the coefficient vectors of the two complexes $B^-$ and $B^+$.
    The two reactions
    $$
        A\xlongrightarrow{\lambda_a\cdot L}B^-,\qquad A\xlongrightarrow{\lambda_a}B^+
    $$
    contribute 
    $$
        \lambda_a(L(b^--a)+(b^+-a))=\lambda_a(La-Lu-La+a+p_a+Lu-a)=\lambda_ap_a=c_a
    $$
    to the coefficient of the monomial $x^a$ in the mass-action ODE and thus realize $c_ax^a$.
    Hence, $k_a=k(a,c_a)$ in all cases.
    All product complexes given in this proof are fixed by $G_a$.
    Thus, taking all $G$-translates of the reactions constructed this way gives a sparsest $G$-symmetric CRN realizing $f$ by Proposition \ref{prop:sparsest}.
\end{proof}

Theorem \ref{thm: sparsest} finds a sparsest CRN only if all coefficient vectors are rational.
This result extends to real coefficients as follows.
\begin{theorem}\label{rmk:real}
Let $f \in \RR[x]^n$ be $G$-equivariant and realizable by a CRN, and let $a \in S/G$.
For the coefficients $c_a=(c_1,\ldots,c_m)$ of a monomial $x^a$ in $f$, let $r=\dim_\QQ\operatorname{span}_\QQ\setn{c_1,\ldots,c_m}$, $D_{c_a}=\setcond{v\in\RR^m}{\forall q\in \QQ^m: q\cdot c_a=0\Rightarrow q\cdot v=0}$ and $V(a,c_a)=\setcond{v\in D_{c_a}\cap\ZZ^m\setminus\setn{0}}{a+v\geq0}$.
Then 
$$
    k(a,c_a)=\begin{cases}
        r, &c_a\in\operatorname{cone}(V(a,c_a)),\\
        r+1, &c_a\notin\operatorname{cone}(V(a,c_a)),
    \end{cases}
$$
where $\operatorname{cone}(A)$ denotes the convex conical hull of $A \subset \RR^m$.
\end{theorem}
\begin{proof}
    Let $c_a=\sum_{i=1}^{k(a,c_a)}\kappa_i(b_i-a)$ for some $b_i\in\NN_0^m\setminus\setn{a}$ and $\kappa_i > 0$.
    Then for every coordinate of $c_a$ we have $(c_a)_j\in\operatorname{span}_\QQ\setcond{\kappa_i}{i\in[k(a,c_a)]}$.
    Thus,
    \begin{equation}\label{eq:step1}
        r\leq\dim_\QQ\operatorname{span}_\QQ\setcond{\kappa_i}{i\in[k(a,c_a)]}\leq k(a,c_a).
    \end{equation}
        
    Let $\varphi_{c_a}:\QQ^m\to\RR,\,q\mapsto q\cdot c_a$.
    Then $(D_{c_a})^\perp=\operatorname{span}_{\RR}\ker(\varphi_{c_a})$.
    Since $\dim_\QQ\operatorname{im}(\varphi_{c_a})=r$ by definition, we get that $\dim_{\QQ}\ker (\varphi_{c_a})=m-r$. Since a $\QQ$-basis of $\ker(\varphi_{c_a})$ is also $\RR$-linearly independent we obtain $\dim_\RR(D_{c_a})=r.$

    Next, assume that $c_a\in\operatorname{cone}(V(a,c_a))$.
    Since $V(a,c_a)\subset D_{c_a}$ and $\dim_\RR (D_{c_a})=r$, Carathéodory's theorem for convex cones gives a sum representation
    $$
        c_a=\sum_{i=1}^r\kappa_iv_i,
    $$
    where $\kappa_i\geq0$ and $v_i\in V(a,c_a)$.
    By setting $b_i=v_i+a$ for all $i$ with $\kappa_i > 0$, we find that $r\geq k(a,c_a)$.

    In the next step, we show that if $r=k(a,c_a)$, then $c_a\in\operatorname{cone}(V(a,c_a))$.
    So assume that $r=k(a,c_a)$, i.e., that there are $\kappa_i > 0$ and $b_i\in\NN_0^m\setminus\setn{a}$ such that $c_a=\sum_{i=1}^r\kappa_i(b_i-a)$ and set $v_i=b_i-a\neq0$.
    Then by \eqref{eq:step1} we have $r=\dim_\QQ\operatorname{span}_\QQ\setcond{\kappa_i}{i\in[r]}$, i.e., the set $\setcond{\kappa_i}{i\in[r]}$ is $\QQ$-linearly independent.
    Let $q\in\ker(\varphi_{c_a})$.
    We have
    $$
        0=q\cdot c_a=\sum_{i=1}^r\kappa_i(q\cdot v_i).
    $$
    As $q\cdot v_i\in\QQ$, linear independence of the $\kappa_i$ yields that $q\cdot v_i=0$, i.e., $v_i\in D_{c_a}$, for all $i$.
    Since $v_i+a=b_i\in\NN_0^m\setminus\setn{a}$, we even have $v_i\in V(a,c_a)$ and consequently $c_a\in\operatorname{cone}(V(a,c_a))$. Together with \eqref{eq:step1} this shows $k(a,c_a) \geq r+1$ if $c_a \not \in \operatorname{cone}(V(a,c_a))$.

    It is now left to show that $k(a,c_a)\leq r+1$ holds in general.
    For this, let $J_0=\setcond{j}{a_j=0}$ and $K=\setcond{w\in\RR^m}{w_j\geq0\text{ for all } j\in J_0}$. Then $(c_a)_j = \sum_i \kappa_i(b_i)_j \geq 0$ for all $j \in J_0$ by the representation of $c_a$ at the beginning of the proof, and so $c_a\in C\coloneqq K\cap D_{c_a}$. Since $D_{c_a}$ is defined by rational linear equations the set $C$ is a rational polyhedral cone and as such generated by finitely many vectors in $C\cap\ZZ^m$. 
    Since $C \subset D_{c_a}$ and $\dim_{\RR}(D_{c_a})=r$ we obtain from Carathéodory's theorem 
    \begin{equation}\label{eq:step2}
        c_a=\sum_{i=1}^r\lambda_iw_i,
    \end{equation}
    where $\lambda_i\geq0$ and $0\neq w_i\in C\cap\ZZ^m$ for all $i$.
    If $a=0$, then $J_0=[m]$, so $K=\RR^m_{\geq0}$.
    Then $w_i\geq0$ and so $w_i\in V(a,c_a)$, and \eqref{eq:step2} implies $k(a,c_a)\leq r$.
    Otherwise, for every $i$ let $N_i\geq1$ be a sufficiently large integer so that $w_i+N_ia\geq -a$ and $w_i+N_ia \neq 0$.
    These integers always exist since $(w_i)_j \geq 0$ for all $j \in J_0$ and $a \neq 0$.
    We obtain
    $$
        c_a=\sum_{i=1}^r\lambda_i(w_i+N_ia)+(\sum_{i=1}^r\lambda_iN_i)(-a).
    $$
    With $\kappa_i=\lambda_i$, $b_i=w_i+(N_i+1)a$, $\kappa_{r+1}=\sum_{i=1}^r\lambda_iN_i$ and $b_{r+1}=0$, this is a linear combination of the desired form $c_a=\sum_{i=1}^{r+1}\kappa_i(b_i-a)$ with $b_i \in \NN_0^m \setminus \{a\}$. So we can conclude $k(a,c_a)\leq r+1$.
\end{proof}

\begin{remark}
    Whereas the proof of Theorem \ref{thm: sparsest} is constructive, this is not the case for Theorem \ref{rmk:real}. 
    In fact, deciding whether two real numbers are $\QQ$-linearly dependent is not possible in general
    in the Blum--Shub--Smale model \cite{zbMATH01096865}. Thus in general $r$ cannot be computed from $c_a$ and since $k(a,c_a)=r$ for $a=0$ by the proof of Theorem \ref{rmk:real}, no general construction of a sparsest CRN can exist in this model. 
    On the other hand, if the $\QQ$-linear relations among the coordinates of $c_a$ are known, then $D_{c_a}$ is explicit and the proof constructs a sparsest realizing CRN. This works because a finite generating set of $\operatorname{cone}(V(a,c_a))$ can be computed.
\end{remark}

\section{Steady states of symmetric CRNs} \label{sec: 3}
In this section we present an algorithm which reduces the search space for steady states of a symmetric CRN $(\S,\C,\R)$ with $n$ species to a finite union of intersections of reflection and coordinate hyperplanes.
The main advantage is that these sets are in general lower dimensional than the real affine space $\RR^{n}$. However, not every steady state of the CRN needs to be contained in this search space. The algorithm relies on an iterative application of Timofte's half-degree principle \cite{Riener12,Timofte}. 

We use the formulation from \cite{Riener12} of the half-degree principle for solutions of symmetric polynomial systems of equations in the non-negative orthant. 
For $a \in \RR^n$ let 
\[n^*(a):=|\setcond{a_j}{1\leq j \leq n, a_j \neq 0}|\] 
denote the number of distinct non-zero coordinates of $a$. Furthermore, for $J \subset [n]$ we write $a_J$ for the projection of $a$ onto the subspace whose coordinates are indexed by $J$, e.g. $a_{\{1,2,n\}}=(a_1,a_2,a_n)$, and for a polynomial $g \in \RR[x]$ we write $\deg_Jg$ for the degree of $g$ considered as a polynomial in the variables $x_j$ for $j \in J$ with coefficients in $\RR[x_i : i \not \in J]$. A polynomial $g \in \RR[x]$ is \emph{symmetric} if it is $S_n$-invariant in the sense of Corollary \ref{cor: F invariant} and $V(g):= \setcond{ a \in \RR^n }{ g(a)=0}$ denotes its real variety.

\begin{theorem}[Timofte's half-degree principle, {\cite[Corollary~1.3]{Riener12}}]\label{thm:half-degree}
     Let $g\in \RR[x]$ be a symmetric polynomial of degree at most $d $. Then the real variety $V(g)$ has non-empty intersection with $\RR_{\geq 0}^n$ if and only if $g(a)=0$ for some $a \in \RR_{\geq 0}^n$ with $n^*(a)\leq \max\{\lfloor \frac{d}{2}\rfloor,2\}$.
\end{theorem}

For CRNs the assumption that the full network is invariant under the action of the symmetric group is in general too strict and we consider Young subgroups of the symmetric group. Let $[n] = J_1 \sqcup \dots \sqcup J_l$ be a set partition with $|J_i|=\lambda_i$. The \emph{Young subgroup} associated with it is 
\[S_\lambda:=S_{J_1}\times \ldots \times S_{J_l} \subset S_n,\]
which is the subgroup of all permutations in $S_n$ which map each block $J_i$ onto itself. The restriction of $S_\lambda$ to $J_i$ is the full symmetric group of $J_i$ and $S_\lambda$ does not interchange different blocks. We write $\lambda=(\lambda_1,\dots,\lambda_l)$ for the sequence of block sizes. For instance, with $J_1=\{1,2\}$ and $J_2=\{3\}$ the polynomial $x_1x_2+x_1^4x_3+x_2^4x_3$ is $S_{(2,1)}$-invariant. Since Young subgroups associated with set partitions of $[n]$ with the same sequence of block sizes are conjugate in $S_n$, the ambiguity in the notation $S_\lambda$ will not create any issues. 
\begin{remark}
If $G \subset S_n$ is any subgroup, the transpositions contained in $G$ generate a Young subgroup of $S_n$ contained in $G$ and it is the largest Young subgroup of $S_n$ contained in $G$. This is the Young subgroup to be used in Theorem \ref{thm:iterative half-degree} and Algorithm \ref{alg:bounding complexity steady states}.  
\end{remark}

\begin{theorem}\label{thm:iterative half-degree}
Let $S_\lambda = S_{J_1}\times \dots \times S_{J_l}$ be a Young subgroup of $S_n$ and let $\N = (\S, \C, \R)$ be a $G$-symmetric CRN with mass-action ODE $(f_1,\ldots,f_n)$, where $G\subset S_n$ is a group containing $S_\lambda$. 
Let $F:=f_1^2+\dots+f_n^2$, $d_i := \deg_{J_i}F$ and $k_i := \max \left\{ \lfloor \frac{d_i}{2} \rfloor,2 \right\}$ for $1 \leq i \leq l$.
Then $\N$ has a steady state if and only if it has a steady state in the set
\[ \mathcal{H} := \left\{ a \in \RR^n_{\geq 0}  \mid n^*\left(a_{J_i }\right) \leq k_i, 1 \leq i \leq l\right\}.\]
\end{theorem}
\begin{proof}
    By Proposition \ref{prop: symmetric CRN implies mass-action ODE equivariant} the mass-action ODE $f$ is $G$-equivariant. So by Corollary \ref{cor: F invariant} the polynomial $F$ is $G$-invariant and thus also $S_\lambda$-invariant. Moreover, $F(a)=0$ if and only if $f(a)=0$. Note that it suffices to prove that if $\mathcal{N}$ has a steady state $p \in \RR_{\geq 0}^n$ then it must also have a steady state in $\mathcal{H}$.
    Without loss of generality we can assume $J_{1} = \{1,\ldots,\lambda_1\}, J_{2} = \{\lambda_1+1,\ldots,\lambda_1+\lambda_2\},\ldots, J_{l} = \{\lambda_1+\ldots+\lambda_{l-1}+1,\ldots,n\}$. 
    
    Since $F$ is $S_\lambda$-invariant and the permutations of $J_l$ fix the remaining coordinates, the polynomial 
    \[g(y):=F(p_1,\ldots,p_{n-\lambda_l},y_1,\ldots,y_{\lambda_l}) \in \RR[y_1,\dots,y_{\lambda_l}]\] is symmetric in $y$. The degree of $g$ is at most $ d_l$ and $g(p_{J_l})=0$ with $p_{J_l} \in \RR_{\geq 0}^{\lambda_l}$. By Theorem \ref{thm:half-degree} there is $\tilde{z} \in \RR_{\geq 0}^{\lambda_l}$ with $g(\tilde{z})=0$ and $n^*(\tilde{z}) \leq k_l$. Then $\tilde{p}=(p_1,\dots,p_{n-\lambda_l},\tilde{z})$ is a steady state with $n^*(\tilde{p}_{J_l}) \leq k_l$. We can proceed by considering the block $J_{l-1}$, fixing all coordinates of $\tilde{p}$ outside $J_{l-1}$, then $J_{l-2}$ and so on. Each step changes only the coordinates of the block under consideration. So after $l$ steps we obtain a steady state in $\mathcal{H}$.
\end{proof}
Moreover, if $P \subset \RR^n$ is an $S_\lambda$-invariant affine subspace then one can apply Theorem \ref{thm:iterative half-degree} also to steady states in $P$. This is because $\operatorname{dist}(x,P)^2$ is a polynomial of degree at most $2$ which is non-negative, $S_\lambda$-invariant and vanishes exactly on $P$. So we can add $\operatorname{dist}(x,P)^2$ to $F=f_1^2+\dots+f_n^2$ and follow the proof of the theorem since the values $k_i$ do not change.

In Section \ref{sec: 4} we present three examples where Theorem \ref{thm:iterative half-degree} is applied.

\begin{remark}\label{rem:molecularity}
Since $F=f_1^2+\dots+f_n^2$ is a sum of squares and $d_i=2\max_{1 \leq j \leq n} \deg_{J_i}f_j$ we have $k_i=\max\{\deg_{J_i}f,2\}$, where $\deg_{J_i}f:=\max_j \deg_{J_i}f_j$. 
This is because the leading forms cannot cancel in a sum of squares.
For a mass-action ODE every monomial of $f$ is of the form $x^a$ with $A \in \mathcal{C}$ a reactant and so $\deg_{J_i}f$ is at most the maximal \emph{$J_i$-molecularity} 
\begin{align}\label{eq: Molecularity}
    M_i := \max\left\{\sum_{j \in J_i} a_j : A \to B \in \R\right\}
\end{align}
of a reactant complex. 
So Theorem \ref{thm:iterative half-degree} says that if $\N$ has a steady state, then it has one in which the species of each block $J_i$ take at most $\max\{M_i,2\}$ distinct non-zero concentrations.
\end{remark}

The set $\mathcal{H}$ is a finite union of intersections of $\RR_{\geq 0}^n$ with linear subspaces which are intersections of coordinate hyperplanes $x_j=0$ and reflection hyperplanes $x_j=x_k$ with $j,k$ in the same block of the set partition of $[n]$. Since $f$ is $G$-equivariant, the set of steady states is $S_\lambda$-invariant and it suffices to search one subspace in each $S_\lambda$-orbit of subspaces in $\mathcal{H}$.

These orbits can be described by patterns. For each block $J_i$ a \emph{pattern} is a tuple of non-negative integers $(m_0;m_1,\dots,m_q)$ with $q \leq k_i$, $m_1 \geq \dots \geq m_q \geq 1$ and $\sum_{j=0}^q m_j = |J_i| = \lambda_i$.
The set corresponding to a pattern consists of the vectors whose block $J_i$ has $m_0$ zero coordinates and $m_j$ coordinates equal to a common value $z_j >0$ for $1 \leq j \leq q$ and the points $z_1,\dots,z_q$ are pairwise distinct. As the representative of its $S_\lambda$-orbit we take the subspace in which the first $m_0$ coordinates of the block are zero, the next $m_1$ coordinates are $z_1$ and so on. For instance, if there is only one block $J_1=[5]$, then the pattern $(2;2,1)$ encodes the set $\setcond{(0,0,a,a,b) \in \RR^5}{a,b >0, a \neq b}$.

For non-negative integers $k,m$ let $\pi_{\leq k}(m)$ denote the number of partitions of $m$ into at most $k$ parts. For fixed $k$ the function $\pi_{\leq k}(m)$ increases polynomially in $m$. 
The number of patterns which need to be considered for the block $J_i$ is $\sum_{m_0=0}^{\lambda_i} \pi_{\leq k_i}(\lambda_i-m_0)$.

Restricting the equation $f=0$ to a pattern for each block gives a reduced system in the unknowns $z_j^{(i)}$ for $1 \leq i \leq l$ and $1 \leq j \leq q^{(i)}$. 
If one solves the reduced systems for all patterns one obtains all steady states that are contained in $\mathcal{H}$ together with their $S_\lambda$-orbits.
This procedure is summarized in Algorithm \ref{alg:bounding complexity steady states} which assumes finitely many solutions in $\mathcal{H}$. 

\begin{center}
   \begin{algorithm} 
\caption{Searching for steady states in a symmetric CRN in the set $\mathcal{H}$ of Theorem \ref{thm:iterative half-degree}}\label{alg:bounding complexity steady states}
\KwData{$G$-symmetric CRN $\N = (\S, \C, \R)$ with mass-action ODE $f=(f_1,\ldots,f_n)$}
\KwResult{All steady states of $\N$ in $\mathcal{H}$ up to the action of $S_\lambda$}
Compute $S_\lambda = S_{J_1} \times \ldots \times S_{J_l}$ the Young subgroup generated by all transpositions in $G$ \\
\For{$1 \leq i \leq l$}{Set $k_i := \max \{\deg_{J_i}f,2\}$} 
Set $L:=\emptyset$\\
\For{each pattern $(m_0^{(i)};m_1^{(i)},\dots,m_{q^{(i)}}^{(i)})$ with $q^{(i)} \leq k_i$ and one for each block $J_i$}{
Substitute the patterns into $f=0$ \\
Compute the real solutions $(z_j^{(i)})$ with $z_j^{(i)} > 0$ and $z_1^{(i)},\dots,z_{q^{(i)}}^{(i)}$ pairwise distinct for $1 \leq i \leq l$ \\
Append each solution $(z_j^{(i)})$ and its pattern to $L$
}
\textbf{return} $L$ 
\end{algorithm} 
\end{center}

Sometimes one is only interested in positive steady states, i.e., solutions of $f=0$ in the positive orthant $\RR_{>0}^n$. 
In this case we cannot apply the half-degree principle.
As Example \ref{example:35} illustrates, the variety $V(g)$ of a symmetric polynomial $g \in \RR[x]$ of degree $d$ can contain points from the positive orthant but the intersection $V(g) \cap \RR_{>0}^n$ need not contain any point with at most $\max\{\lfloor \frac{d}{2} \rfloor,2\}$ distinct coordinates.

\begin{example}\label{example:35}
    Consider the symmetric polynomial $g(x_1, x_2, x_3) := e_1^4 - 15e_1^2 e_2 + 50e_2^2 + 50e_1 e_3,$ where $e_1= x_1+x_2+x_3, e_2 = x_1x_2+x_1x_3+x_2x_3, e_3 = x_1x_2x_3 \in \RR[x_1,x_2,x_3]$.
By Timofte's half-degree principle \cite{Riener12,Timofte} $g$ has a zero in $\RR_{\geq 0}^3$ if and only if it has a root in $\RR_{\geq 0}^3$ with at most two distinct non-zero coordinates. One such root is $\left(1+\sqrt{\frac{3}{5}},1-\sqrt{\frac{3}{5}},0\right)$. 
Moreover, $g$ has zeros in the open orthant $\RR_{>0}^3$, e.g., the point whose coordinates are the three distinct roots of $t^3-t^2+\frac{3}{20}t-\frac{1}{400}$ which are real and positive. Then $e_1=1, e_2=\frac{3}{20}$ and $e_3=\frac{1}{400}$ and hence $g = 1-\frac{45}{20}+\frac{450}{400}+\frac{50}{400}=0$. 
However, $g$ has no zero in $\RR_{>0}^3$ with at most two distinct positive coordinates. This can be seen as follows. If $g(a,a,b) = 0$ for some $a, b > 0$ we can assume that $e_1(a,a,b)=1$, since $g$ is homogeneous. We substitute $b=1-2a$ and a computation shows that the roots of $g(a,a,1-2a)=450a^4-700a^3+295a^2-30a+1$ are $\approx 0.5105$ and $\approx 0.9172$ and a pair of complex conjugate roots. However, for $a > 0.5$ we have $1-2a < 0$.
\end{example}

In \cite{riener2018real} the existence of points with few distinct coordinates in non-empty basic semialgebraic sets defined by symmetric equations and equivariant inequalities is proved. The authors use Timofte's half-degree principle for polynomial equations and, for inequalities, a result of Shvartsman \cite{shvartsman1982some} on the structure of equivariant maps in $\RR[x]^n$ with respect to the action of a real reflection group. We can adapt \cite[Corollary~3.4]{riener2018real} and \cite[Lemma~3.7]{riener2018real} to Young subgroups. 

\begin{lemma}\label{lem: Shvartsman}
Let $S_\lambda = S_{J_1} \times \cdots \times S_{J_l} \subset S_n$ be a Young subgroup and let $f \in \RR[x]^n$ be $S_\lambda$-equivariant. For $1 \leq i \leq l$ let $D_i:=\max_{j \in J_i} \deg_{J_i} f_j$ and $p_{J_i,k} := \sum_{j \in J_i}x_j^k$ for $k \in \NN$. 
Then there exist unique $S_\lambda$-invariant polynomials $g_1^{(i)},\dots,g_{|J_i|}^{(i)}$ with $\deg_{J_i} g_k^{(i)} \leq D_i-k+1$ such that for every $j \in J_i$
\begin{align}\label{eq: decomposition}
    f_j(x) = \sum_{k=1}^{|J_i|} g_k^{(i)}(x)\cdot  \frac{\partial p_{J_i,k}}{\partial x_j} = \sum_{k=1}^{|J_i|} k g_k^{(i)}(x) x_j^{k-1}.
\end{align}
If $f( \tilde{x})=0$ then every coordinate $\tilde{x}_j$ with $j \in J_i$ is a root of 
\begin{align}\label{eq: Psi}
    \Psi_{i,\tilde{x}}(u):=\sum_{k=1}^{|J_i|} k g_k^{(i)}(\tilde{x}) u^{k-1}.
\end{align} 
Furthermore, if $\Psi_{i,\tilde{x}} \neq 0$ then $\tilde{x}_{J_i}$ has at most $\min \{D_i,|J_i|-1\}$ distinct coordinates.
\end{lemma}
\begin{proof}
Since the group $S_\lambda$ is a finite reflection group its invariant ring is generated by the algebraically independent power sums $p_{J_i,k_r^{(i)}}$ for $1 \leq i \leq l$ and $1 \leq k_r^{(i)} \leq |J_i|$. 
By a result of Shvartsman \cite{shvartsman1982some} (see also \cite[Theorem~3.3]{riener2018real}) the $S_\lambda$-equivariant function $f$ is a unique combination of the gradients of these power sums with $S_\lambda$-invariant polynomials as coefficients. 
For $1 \leq i,r \leq l$ and $j \in J_i$ we have $\frac{\partial p_{J_r,k}}{\partial x_j} \neq 0$ only if $J_r=J_i$. From the Vandermonde determinant it can be seen that the vectors $(x_j^{k-1})_{j \in J_i}$ for $1 \leq k \leq |J_i|$ are linearly independent over the quotient field of $\RR[x]$.
So we obtain that $\deg_{J_i} g_k^{(i)} \leq D_i-k+1$ which proves \eqref{eq: decomposition}. 

Finally, from $\deg_{J_i}f_j \leq D_i$ we obtain the claim on the upper bound of distinct coordinates of $\tilde{x}_{J_i}$.
\end{proof}
For $l=1$ the decomposition in \eqref{eq: decomposition} coincides with \cite[Corollary~3.4]{riener2018real} and the last statement corresponds to \cite[Lemma~3.7]{riener2018real}.

\begin{corollary}\label{cor:schvartsman}
Let $f \in \RR[x]^n$ be the mass-action ODE of an $S_\lambda=S_{J_1}\times \cdots \times S_{J_l}$-symmetric CRN and let $D_i:=\max_{j \in J_i} \deg_{J_i} f_j$. For each block $1 \leq i \leq l$ let $\Psi_{i,\tilde{x}}(u)$ denote the univariate polynomial in \eqref{eq: Psi}.
If $\tilde{x}$ is a steady state and $\Psi_{i,\tilde{x}} \neq 0$ then the coordinates of $\tilde{x}_{J_i}$ take at most $D_i$ distinct values.
\end{corollary}
\begin{proof}
This is a direct consequence of Lemma \ref{lem: Shvartsman}.
\end{proof}

The corollary actually applies to all real solutions of $f=0$ and not only to steady states of the CRN.
\begin{remark}
Corollary \ref{cor:schvartsman} and Theorem \ref{thm:iterative half-degree} complement each other. We have $D_i \leq \deg_{J_i} f \leq M_i$, where $M_i$ denotes the maximal $J_i$-molecularity defined in \eqref{eq: Molecularity}, but the corollary can only be applied in blocks with $\Psi_{i,\tilde{x}} \neq 0$. A priori we do not know if $\Psi_{i,\tilde{x}}\neq 0$ holds for a steady state but Theorem \ref{thm:iterative half-degree} can still be used. However, Theorem \ref{thm:iterative half-degree} does not give information about steady states outside of the set $\mathcal{H}$. If indeed $\Psi_{i,\tilde{x}} \neq 0$ for all steady states and all $1 \leq i \leq l$ then all steady states of the symmetric CRN are contained in the set $\mathcal{H}$.
\end{remark}

\begin{remark}
One might hope that the use of Specht ideals (see \cite{Moustrou}), which are certain symmetric ideals associated with the irreducible representations of the symmetric group, yields a more refined bound on the combinatorial structure of steady states. 
If a CRN is symmetric under the full symmetric group $S_n$, this approach is effective. Writing a steady state up to permutation in the form 
\[(0, \ldots, 0, a_1, \ldots, a_1, a_2, \ldots, a_2, \ldots, a_l, \ldots, a_l)\]
with pairwise distinct $a_i > 0$, both the number $l$ of distinct positive coordinates and their multiplicities are restricted to a finite set of possibilities characterized by the varieties of Specht ideals (see \cite[Proposition~3]{Moustrou}). This is since the ideal generated by a family of equivariant polynomials, which is the case for symmetric CRNs by Proposition \ref{prop: symmetric CRN implies mass-action ODE equivariant}, is closed under the group action. 

However, extending this framework to Young subgroups presents significant algebraic challenges. When applying these methods iteratively across disjoint groups of variables, one must work with polynomials in a subset of variables over a polynomial ring in the remaining variables. Since this base ring is not a field, the standard theory of Specht ideals does not apply. Consequently, we are not aware of a general method to circumvent the study of syzygies that depend on the equivariant mass-action ODE.
\end{remark}

\subsection{Examples}\label{sec: 4}
We apply Theorem \ref{thm:iterative half-degree} and Algorithm \ref{alg:bounding complexity steady states} to three families of CRNs.
These examples are motivated by the ``origin of life'' literature.

\subsubsection{Bimolecular autocatalytic system}
A \emph{bimolecular autocatalytic system} (see \cite[Definition~3.4]{craciun2022autocatalytic}) is a CRN with species $X_1,\dots,X_n$, where all reactions are of the form $X_i+X_j \to X_i+X_j+X_l$ for some $i,j,l \in \{1,\dots,n\}$.

For $n \geq 3$ we consider the following instance. The reactions are
\begin{align*}
    2X_i &\xlongrightarrow{\kappa_1} 3X_i \quad & 1 \leq i \leq n \\
    X_i + X_j &\xlongrightarrow{\kappa_2} X_i + X_j + X_l \quad &1 \leq i \leq j \leq n,\  l \notin \{i,j\}.
\end{align*}
So the CRN is $S_n$-symmetric.
We write $e_1(x) = \sum_{i=1}^n x_i$ and $p_2(x)= \sum_{i=1}^n x_i^2$.
Let $\mathrm{e}_i$ denote the $i$-th standard basis vector in $\RR^n$.
Since $\sum_{1\leq i\leq j \leq n, \, l \notin \{i,j\}} x_ix_j = \frac{1}{2}( (e_1(x)-x_l)^2 + p_2(x) -x_l^2 )$, the mass-action ODE is
\begin{align*}
    f(x) & = \sum_{i=1}^n \kappa_1 x_i^2 \mathrm{e}_i + \sum_{1 \leq i \leq j \leq n, \, l \not\in \{i,j\}} \kappa_2 x_ix_j \mathrm{e}_l  =  \left(
\kappa_1 x_l^2 + \frac{\kappa_2}{2}\left( (e_1(x)-x_l)^2 + p_2(x) - x_l^2 \right)\right)_{1 \leq l \leq n}.
\end{align*}
Observe that each coordinate function $f_l$ is positive on $\RR_{\geq 0}^n \setminus\{0\}$ and so the only steady state is $0$.
It is natural to study the relative populations $\frac{1}{e_1(x)}x$ as $e_1(x)$ is the total concentration of the mass. 
By \cite[Theorem~3.5]{craciun2022autocatalytic}, after the time rescaling
$\mathrm{d}\tau = e_1(x(t))\,\mathrm{d}t$, the relative populations $\frac{1}{e_1(x)}x$
evolve according to the mass-action system of the CRN obtained by adding one copy of each species $X_p$ with $p \neq l$ to the reactant complex of every reaction $X_i+X_j \to X_i+X_j+X_l$.
The new CRN is
\begin{align*}
    X_p+2X_i &\xlongrightarrow{\kappa_1} 3X_i \quad & 1 \leq i \leq n, \ p \neq i \\
    X_p+X_i + X_j &\xlongrightarrow{\kappa_2} X_i + X_j + X_l \quad &1 \leq i \leq j \leq n,\  l \notin \{i,j\}, \ p \neq l
\end{align*}
which is again $S_n$-symmetric, so we have one invariant block of variables. By \cite[Corollary~3.3]{craciun2022autocatalytic} its mass-action ODE $g$ is
\[g(x) = e_1(x) f(x) - \sum_{i=1}^n f_i(x) x. \]
Since $g$ is homogeneous and $\sum_{i=1}^n g_i = 0$, the total concentration $e_1(x)$ of the new CRN is conserved and we may restrict to solutions in the simplex $e_1(x)=1$.

We apply Theorem \ref{thm:iterative half-degree} and Remark \ref{rem:molecularity}. The maximal molecularity of the new CRN is $3$. Thus there exists a steady state if and only if there exists one with $n^*(x) \leq \max\{3,2\}=3$ distinct non-zero coordinates.

In the simplex $e_1(x)=1$ we have
\[ g_l = \kappa_1 x_l^2 - (\kappa_2+\lambda)x_l + \frac{\kappa_2}{2}\left( 1 + p_2(x) \right), \text{ where } \lambda := \sum_{i=1}^n f_i(x) = \kappa_1 p_2(x) + \frac{\kappa_2}{2}\left( n-2 + n\,p_2(x) \right). \]
In particular, every coordinate of a steady state $\tilde{x}$ is a root of the same univariate quadratic polynomial
\[ \Psi(u) := \kappa_1 u^2 - (\kappa_2+\lambda)u + \frac{\kappa_2}{2}\left( 1 + p_2(\tilde{x}) \right). \]
Since also $\Psi(0)>0$, every steady state lies in $\RR_{>0}^n$ and has at most two distinct coordinates. Therefore it lies in the set $\mathcal{H}$ and Algorithm \ref{alg:bounding complexity steady states} finds all steady states.
The only patterns that can have solutions are $(0;n)$ and $(0;m_1,m_2)$ with $m_1 \geq m_2 \geq 1$ and $m_1+m_2=n$. For instance, the pattern $(0;n)$ gives the steady state $\left( \frac{1}{n},\dots,\frac{1}{n} \right)$ for all $\kappa_1,\kappa_2>0$.

\subsubsection{Bipartite reversible autocatalytic network}
In \cite{Deshpande,Gopalkrishnan} the autocatalytic CRN 
\[ 
Y \longleftrightarrow 2X, \qquad X \longleftrightarrow 2Y
\]
is considered.
We generalize this CRN to a bipartite network consisting of two blocks of species $\mathcal{X} = \{X_1, \dots, X_n\}$ and $\mathcal{Y} = \{Y_1, \dots, Y_m\}$ with reactions
\begin{align*}
    Y_j &\xrightleftharpoons[\kappa_2]{\kappa_1} 2X_i \quad & 1 \leq i \leq n, 1 \leq j \leq m \\
    2Y_j &\xrightleftharpoons[\kappa_4]{\kappa_3} X_i \quad & 1 \leq i \leq n, 1 \leq j \leq m
\end{align*}
where $n,m \geq 1$.
The CRN is $ S_n \times S_m$-symmetric, where $S_n$ acts on $\mathcal{X}$ and $S_m$ acts on $\mathcal{Y}$. We write $x_i$ for the concentration of $X_i$, $y_j$ for the concentration of $Y_j$ and $e_1(y)=\sum_{j=1}^m y_j, p_2(y)=\sum_{j=1}^m y_j^2$ and analogously for $x$. The CRN has the mass-action ODE $f=(f_{X_1},\dots,f_{X_n},f_{Y_1},\dots,f_{Y_m})$ with
\begin{align*}
  f_{X_i}(x,y) &= \sum_{j=1}^m \left( 2\kappa_1 y_j - 2\kappa_2 x_i^2 + \kappa_3 y_j^2 - \kappa_4 x_i \right) = -2m\kappa_2x_i^2-m\kappa_4x_i + 2\kappa_1 e_1(y)+\kappa_3p_2(y) \\   f_{Y_j}(x,y) &= \sum_{i=1}^n \left( -\kappa_1 y_j + \kappa_2 x_i^2 - 2\kappa_3 y_j^2 + 2\kappa_4 x_i \right) = -2n\kappa_3y_j^2-n\kappa_1y_j + 2\kappa_4 e_1(x)+\kappa_2p_2(x) .
\end{align*} 
The variables partition into the blocks $J_1 = \{x_1, \dots, x_n\}$ and $J_2 = \{y_1, \dots, y_m\}$. 
The maximal $J_1$- and $J_2$-molecularities are both $2$. So by Theorem \ref{thm:iterative half-degree} and Remark \ref{rem:molecularity} the system has a steady state if and only if there exists a steady state where the number of distinct non-zero coordinates within each block satisfies:
\[ n^*(x) \leq \max \left\{ 2, 2 \right\} = 2, \qquad n^*(y) \leq \max \left\{ 2, 2 \right\} = 2. \]
In this example $n^*(x) \leq 1$ and $n^*(y) \leq 1$ hold for every steady state. At a steady state $(\tilde{x},\tilde{y})$ all $\tilde{x}_i$ are roots of the same quadratic polynomial
\[ \Psi(u) := -2m\kappa_2 u^2- m \kappa_4u+2\kappa_1e_1(\tilde{y})+\kappa_3p_2(\tilde{y}).\]
This polynomial is strictly decreasing on $\RR_{\geq 0}$ and so it has exactly one non-negative root. This root is $0$ if and only if $\tilde{y}=0$ and otherwise $\left( \frac{\kappa_1^2 \kappa_4}{\kappa_2^2 \kappa_3} \right)^{\frac{1}{3}}$. 
A similar argument shows that $\tilde{y}=0$ if and only if $\tilde{x}=0$ and otherwise $\tilde{y}_j = \left( \frac{\kappa_1 \kappa_4^2}{\kappa_2 \kappa_3^2} \right)^{\frac{1}{3}}$ for $1 \leq j \leq m$.
Consequently, the only patterns with solutions are $(n)$ for $J_1$ together with $(m)$ for $J_2$ which corresponds to the steady state $0$, and the patterns $(0;n)$ and $(0;m)$ which give the positive steady state
$$x_i = \left( \frac{\kappa_1^2 \kappa_4}{\kappa_2^2 \kappa_3} \right)^{\frac{1}{3}} \quad 1 \leq i \leq n, \qquad y_j = \left( \frac{\kappa_1 \kappa_4^2}{\kappa_2 \kappa_3^2} \right)^{\frac{1}{3}} \quad 1 \leq j \leq m.$$

\subsubsection{Mutual inhibition and homodimer self-activation}\label{subsec:mutual}

\begin{figure}
	\begin{center}
		\includegraphics{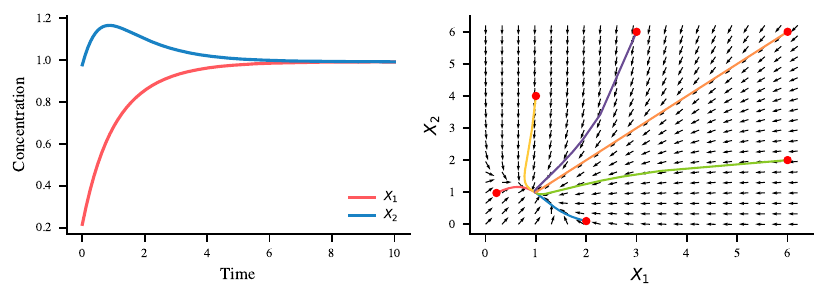}
	\end{center}
	\caption{Simulation of Example \ref{subsec:mutual} for $n=2$ with randomly sampled parameters $\kappa_i$, $1 \leq i \leq 5$. The right panel shows the phase portrait of the corresponding ODE system, including six solution trajectories. The red dots indicate the initial conditions of these trajectories.}
	\label{fig:example43}
\end{figure}

We consider $n$ copies of the Schlögl model \cite{schlogl1972chemical}, the classical one-species mass-action model of bistability, coupled by pairwise annihilation:
\begin{align*}
    \emptyset &\xrightleftharpoons[\kappa_2]{\kappa_1} X_i & \qquad 1 \leq i \leq n \\
    2X_i &\xrightarrow{\kappa_3} 3X_i & \qquad 1 \leq i \leq n \\
    3X_i &\xrightarrow{\kappa_4} 2X_i & \qquad 1 \leq i \leq n \\
    X_i+X_j & \xrightarrow{\kappa_5} \emptyset & \qquad 1 \leq i < j \leq n.
\end{align*}
These reactions represent basal influx and first-order decay, cooperative autocatalysis with trimolecular saturation and mutual inhibition of competing autocatalytic species.
Without the last reaction the mass-action ODE would consist of $n$ copies of the same cubic equation. The annihilation reaction preserves the $S_n$-symmetry but couples the variables. So we have just one invariant block of variables.
The mass-action ODE is
\[f=(f_1,\dots,f_n)=\left(\kappa_1 - \kappa_2 x_i + \kappa_3 x_i^2 - \kappa_4 x_i^3 - \kappa_5 x_i \sum_{j \neq i} x_j\right)_{1 \leq i \leq n}.\]
 Figure~\ref{fig:example43} shows one simulation together with the corresponding phase portrait.
 
We apply Theorem \ref{thm:iterative half-degree} and Remark \ref{rem:molecularity}. The maximal molecularity is $3$.
Thus the CRN has a steady state if and only if there exists one with $n^*(x) \leq \max\{3,2\}=3$ distinct non-zero coordinates.
Since $\kappa_1 > 0$ every steady state must lie in $\RR_{>0}^n$. 
Setting $s = \sum_{j=1}^n x_j$ we observe
\[ f_i = \kappa_1-(\kappa_2+\kappa_5s)x_i+(\kappa_3+\kappa_5)x_i^2-\kappa_4x_i^3.\]
In particular, every coordinate of a steady state is a root of the same univariate cubic polynomial
\[ \Psi_s(u):=\kappa_1-(\kappa_2+\kappa_5s)u+(\kappa_3+\kappa_5)u^2-\kappa_4u^3,\]
where $s$ is the sum of the coordinates of the steady state.
So every steady state has at most three distinct coordinates and therefore lies in the set $\mathcal{H}$. Thus, Algorithm \ref{alg:bounding complexity steady states} finds all steady states. 
The only patterns that can have solutions are $(0;n), (0;m_1,m_2)$ with $m_1 \geq m_2 \geq 1$ and $m_1+m_2=n$ and 
$(0;m_1,m_2,m_3)$ with $m_1\geq m_2\geq m_3\geq 1$ and $m_1+m_2+m_3=n$. These correspond to partitions of $n$ into at most three parts. 

For the pattern $(0;n)$ the reduced system is the single cubic equation $\Psi_{nz}(z)=0$. For every other pattern the coordinates of steady states are roots of $\Psi_s$. By Vieta's formulas the reduced system is equivalent to
\[ z_1+z_2+z_3=\frac{\kappa_3+\kappa_5}{\kappa_4},\qquad z_1z_2z_3=\frac{\kappa_1}{\kappa_4},\qquad \kappa_4(z_1z_2+z_1z_3+z_2z_3)-\kappa_5(m_1z_1+m_2z_2+m_3z_3)=\kappa_2, \]
where only the last equation depends on the pattern. 
For instance, with $(\kappa_1,\dots,\kappa_5)=(6,5,5,1,1)$ and $s=6$ we have $\Psi_6(u)=-(u-1)(u-2)(u-3)$. For $n=3$ the point $(1,2,3)$ is a steady state with three distinct coordinates. For $n= 6$ the only the steady state is $(1,\dots,1)$ while for $n \geq 7$ there is none (for the fixed concentration $s$ and reaction parameters).

\subsection*{AI Declaration}
This project was initiated in 2023, prior to the widespread adoption of large language models in mathematical research. During the final preparation of the manuscript, the authors utilized large language models (Claude 5.1 Fable, ChatGPT 5.6 Luna) for streamlining proofs and improving the presentation of the paper. An initial proof of Theorem \ref{rmk:real} was suggested by ChatGPT 5.6 Luna.
The final manuscript was written entirely by the authors who take full responsibility for the content.

\subsection*{Acknowledgment}
We thank Alicia Dickenstein for bringing our attention to autocatalytic recombination systems as CRNs with inherent symmetries and Marta Panizzut for suggesting to investigate the application of Specht ideals in the study of steady states of CRNs.\\

This work was supported in part by the Croatian Science Foundation under the project number HRZZ-IP-2024-05-2130.

\printbibliography
\end{document}